\documentclass{amsart}
\newcommand{\note}{\noindent {\bf Notation. }}
\newcommand{\remark}{\noindent {\bf Remark. }}

\newcommand{\ws}{\hspace{4pt}}
\newtheorem{theorem}{Theorem}

\newtheorem{proposition}{Proposition}

\newtheorem{lemma}{Lemma}
\newtheorem{defi}{Definition}
\usepackage[]{amssymb,amsopn}
\usepackage[]{amsmath}
\usepackage[]{eucal}
\usepackage[]{latexsym}

\begin{document}

\title[]{Discrete diffusion semigroups associated with Dunkl-Jacobi and exceptional Jacobi polynomials}
\author{\'A. P. Horv\'ath }

\thanks{Supported by the NKFIH-OTKA Grants K128922 and K132097.}

\begin{abstract}
Some weighted inequalities for the maximal operator with respect to the discrete diffusion semigroups associated with exceptional Jacobi and Dunkl-Jacobi polynomials are given. This setup allows to extend the corresponding results obtained for discrete heat semigroup recently to richer class of differential-difference operators.

\end{abstract}
\maketitle

\section{Introduction}

Diffusion semigroups in the Stein's sense (cf. \cite{ste}) were investigated by several authors. Just some recent examples are mentioned here. For instance in \cite{ab} the dynamics of the heat semigroup generated by the Jacobi operator is studied, in \cite{cgm}  via examination of the (continuous) heat equation associated with the Jacobi-Dunkl operator on the real line, a Poisson equation is solved and a new family of one-dimensional Markov processes is introduced.\\
Besides the continuous operator semigroup, lately the study of discrete diffusion semigroups has come to the forefront of interest, see e.g. \cite{cgrtv}, \cite{abr}, \cite{acl}, \cite{bcf}.

Below we investigate discrete diffusion semigroups related to the recurrence relations of Dunkl-Jacobi and exceptional Jacobi polynomials. The derivation of these two orthogonal systems from the classical Jacobi polynomials shows some similarity. Indeed, Dunkl-Jacobi polynomials are eigenfunctions of a differential-difference operator (see \eqref{djo}) and exceptional Jacobi polynomials are eigenfunctions of a differential operator (see \eqref{ka}), such that the original differential operator of classical Jacobi polynomials can be expressed by the ones mentioned above, see \eqref{DJ} and \eqref{si}, respectively. These systems are complete, but counter to the standard orthogonal polynomials, they possess recurrence formulae with more than three terms. This property allows to extend the examination from the standard discrete heat semigroup to more complicated differential-difference operators. The norm estimates given for the maximal operator of the semigroups show the dependence of the norm of the solution on the norm of the initial value.

This paper is organized as follows: in the next section we introduce the discrete diffusion semigroup, and make some remarks about its general properties. In the third section, by a general theorem proved in \cite{bcf}, we extend the maximal operator of the diffusion semigroup associated with exceptional Jacobi polynomials to weighted $l^p$ spaces. In the last section we deduce a similar theorem for discrete Dunkl-Jacobi semigroup from the corresponding result of \cite{acl}.

\section{The discrete diffusion semigroup}

\subsection{Elementary setup}
Let $I$ be a bounded real interval, and $\mu$ a positive measure supported on $I$ with infinitely many points in its support. Let $\{\phi_n\}_{n\in \mathcal{I}}$ be a (complex valued) complete orthonormal system in $L^2_\mu$, where $\mathcal{I}=\mathbb{N}$ or  $\mathcal{I}=\mathbb{Z}$.
Let $s: I \to \mathbb{R}$ be in $L^\infty_\mu$ on $I$ and  $s_M:=\rm{essup}_Is$. Let us define a kernel function as
\begin{equation}\label{K} K_t(n,m):=\int_Ie^{-(s_M-s(x))t}\phi_n(x)\overline{\phi}_m(x)d\mu(x),\end{equation}
where $t\ge 0$. Let $f\in l^2(\mathbb{C})$ that is a complex valued sequence in $l^2$ and define the operator acting on $f$ as
\begin{equation}W_tf(n):=\sum_{m\in \mathcal{I}}f(m)K_t(n,m).\end{equation}
Since $|e^{-(s_M-s(x))t}|\le 1$, by orthonormality
$$|K_t(n,m)|\le 1.$$
Moreover, considering that
\begin{equation}\label{c} K_t(n,m)=c_m(e^{-(s_M-s)t}\phi_n),\end{equation}
where $c_m(F)$ means the $m^{th}$ Fourier coefficient of $F$ with respect to $\{\phi_n\}$, Parseval's formula gives that
\begin{equation}\label{2n}\|K_t(n,\cdot)\|_2 = \|e^{-(s_M-s)t}\phi_n\|_{2,\mu}\le 1.\end{equation}
This implies that
$$|W_tf(n)|\le \|f\|_2 \|K_t(n,\cdot)\|_2.$$
Thus $W_tf(n)$ is well-defined for each sequence in $l^2$. Let us define the maximal operator
\begin{equation}\label{max} W_*f(n):= \sup_{t>0}|W_tf(n)|,  \ws \ws n\in \mathcal{I}.\end{equation}
First we observe that the family of operators $\{W_t\}_{t\ge 0}$ is a strongly continuous semigroup on $l^2$ and possesses the contraction property wich implies the boundedness of the maximal operator.

\medskip

\begin{proposition} Let $f \in l^2$.
\begin{equation}\label{2}W_0f(n)=f(n) \ws\ws \ws \forall n\in \mathcal{I}.\end{equation}
\begin{equation}\label{1}\|W_tf\|_2\le \|f\|_2.\end{equation}
\begin{equation}\label{3}\sum_{j\in  \mathcal{I}}K_{t_1}(n,j)\overline{K}_{t_2}(m,j)=K_{t_1+t_2}(n,m).\end{equation}
\begin{equation}\label{4}W_{t_1}W_{t_2}f(n)=W_{t_1+t_2}f(n) \ws \ws \ws \forall t_1, t_2\ge 0, \ws n\in \mathcal{I}.\end{equation}
\begin{equation}\label{5}\lim_{t\to 0+}\|W_tf-f\|_2=0.\end{equation}
\begin{equation}\label{7} \lim_{t\to 0+}W_tf(n)=f(n) \ws\ws n\in \mathcal{I}.\end{equation}
\begin{equation}\label{6}\|W_*f\|_2 \le C\|f\|_2,\end{equation}
where $C$ is an absolute constant.
\end{proposition}

\medskip

The properties listed above are well-known, for sake of completeness we prove it in brief.

\proof
Orthogonality implies \eqref{2}.\\
Let us denote by $P$ the set of sequences with finitely many nonzero elements. By density it is enough to prove \eqref{1} for $f\in P$. Let $f(n)=0$, if $|n|>N$, say. Then
$$\|W_t(f)\|_2=\left(\sum_{n\in \mathcal{I}}\left|\sum_{|m|\le N}f(m)K_t(n,m)\right|^2\right)^{\frac{1}{2}}=(*).$$
Denoting by $p_N:=\sum_{|m|\le N}f(m)\phi_m$ and by $p_N^*:=\sum_{|m|\le N}\overline{f}(m)\phi_m$
$$(*)=\left(\sum_{n\in \mathcal{I}}\left|c_n(e^{-(s_M-s)t}p_N^*)\right|^2\right)^{\frac{1}{2}}=\|e^{-(s_M-s)t}p_N^*\|_{2,\mu}\le\|p_N^*\|_{2,\mu}=\|\overline{f}\|_2=\|f\|_2.$$
In view of \eqref{c}
$$\sum_{j\in  \mathcal{I}}K_{t_1}(n,j)\overline{K}_{t_2}(m,j)=\sum_{j\in  \mathcal{I}}c_j(e^{-(s_M-s)t_1}\phi_n)\overline{c}_j(e^{-(s_M-s)t_2}\phi_m)$$ $$=\int_Ie^{-(s_M-s)(t_1+t_2)}\phi_n\overline{\phi}_md\mu=K_{t_1+t_2}(n,m).$$
\eqref{4} is ensured by \eqref{3}:
$$W_{t_1}W_{t_2}f(n)=\sum_{m\in \mathcal{I}}\sum_{j\in \mathcal{I}}f(j)K_{t_2}(m,j)K_{t_1}(n,m)=\sum_{j\in \mathcal{I}}f(j)\sum_{m\in \mathcal{I}}\overline{K}_{t_2}(j,m)K_{t_1}(n,m)$$ $$=\sum_{j\in \mathcal{I}}f(j)K_{t_1+t_2}(n,j),$$
where the last but one equality fulfils because \eqref{1} implies that for each subsequence $\{j_l\}$ the series is convergent.\\
Again, it is enough to prove \eqref{5} for $f\in P$. Let $p_N$ be as above.
$$\|W_tf-f\|_2^2=\sum_{n\in \mathcal{I}}\left|\int_I\left(e^{-(s_M-s)t}-1\right)p_N^*\overline{\phi}_nd\mu\right|^2=\left\|\left(e^{-(s_M-s)t}-1\right)p_N^*\right\|_{2,\mu}^2$$ $$\le \left\|\left(e^{-(s_M-s)t}-1\right)\right\|_\infty^2\|f\|_{2,\mu}^2.$$
Taking the limit, \eqref{5} is proved. \\
\eqref{7} is a consequence of \eqref{5}.\\
According to \cite[Ch. III, p. 73]{ste} \eqref{6} follows from \eqref{1}.

\medskip

\remark
Since $\overline{K}_t(n,m)=K_t(m,n)$, \eqref{4} implies that the operator $W_t$ is positive (definite) for all $t\ge 0$, indeed, $W_t=W_{\frac{t}{2}}^2$.

In accordance with \cite[p. 65]{ste} a family of selfadjoint operators with properties \eqref{2}, \eqref{4}, \eqref{1} is called a symmetric diffusion semigroup.

\subsection{Diffusion semigroup generated by recurrence formulae}

Consider the $L^2_\mu$ space on a finite interval $I$ with the complete orthonormal system  $\{\phi_n\}_{n\in \mathcal{I}}$ as previously. For a function $s\in L^\infty_\mu$ we introduce a multiplication operator
\begin{equation}M_s: L^2_\mu \to L^2_\mu; \ws \ws M_s F:=sF.\end{equation}
Introduce the operator
\begin{equation}A:=s_MI-M_s.\end{equation}
Describing the operators above in the Schauder basis $\{\phi_n\}_{n\in \mathcal{I}}$ one can consider $M_s$,$A$ as operators acting on $l^2$ as well.

By spectral theorem the $W_t$ operators defined above can be expressed as
\begin{equation} W_tf=e^{tA}f=\int_{s(I)}e^{-(s_M-s(\lambda))t}dE_{M_s}(\lambda)f.\end{equation}

We investigate the next initial-value problem.
\begin{equation}\label{h}\left\{\begin{array}{ll}\frac{\partial u(n,t)}{\partial t}=Au(n,t),\\ u(n,0)=f(n),\end{array}\right.\end{equation}
where $A$ is the infinitesimal genarator of the semigroup.

 Now we apply the left-hand side of \eqref{h} to $ W_tf(n)$. Considering \eqref{2n} the derivation can be moved inside, that is
 $$\frac{\partial W_tf(n)}{\partial t}=\sum_{m\in \mathcal{I}}f(m)\int_Ie^{-(s_M-s(x))t}(s(x)-s_M)\phi_n(x)\overline{\phi}_m(x)d\mu(x).$$
 Let $(s(x)-s_M)\phi_n=\sum_{k\in \mathcal{I}}c_{k,n}\phi_k$. That is $c_{k,n}= c_k((s(x)-s_M)\phi_n)$ and by the assumption on $s$, $\{c_{k,n}\}_k$ is in $l^2$. Let $b_{k,m}:=b_k(e^{-(s_M-s)t}\phi_m)$. Then $\{b_{k,n}\}_k$ is in $l^2$ again and by completeness
 $$\int_Ie^{-(s_M-s(x))t}(s(x)-s_M)\phi_n(x)\overline{\phi}_m(x)d\mu(x)=\sum_{k\in \mathcal{I}}c_{k,n}\overline{b}_{k,m}$$ $$=\sum_{k\in \mathcal{I}}c_{k,n}\int_Ie^{-(s_M-s(x))t}\overline{\phi}_m(x)\phi_k(x)d\mu(x).$$
Thus
$$\frac{\partial W_tf(n)}{\partial t}=\sum_{m\in \mathcal{I}}f(m)\sum_{k\in \mathcal{I}}c_{k,n}K_t(k,m).$$
Again by \eqref{1} and Antosik's theorem the order of the sums above can be interchanged. Thus taking into account \eqref{2}, $u(n,t)=W_tf(n)$ solves the initial-value problem \eqref{h}.
Indeed,
$$\frac{\partial W_tf(n)}{\partial t}=\sum_{k\in \mathcal{I}}c_{k,n}W_tf(k)=A W_tf(n).$$

Subsequently we investigate \eqref{h} with special right-hand side.

Let us suppose now that $M_s$ generates a recurrence formula with respect to $\{\phi_n\}_n$ that is $s\phi_n$ can be expressed as a linear combination of certain $\phi_k$-s of constant length for each $n\in \mathcal{I}$. In this case \eqref{h} is an initial-value problem with respect to a discrete partial differential equation, that is on the right-hand side of \eqref{h} there is a discrete differential operator.

For instance, with $u(n,t)=e^{-2t}I_{n-m}(2t)$, where $I_k(t)$ is the Bessel function of imaginary argument, on the right-hand side there is just the discrete Laplacian, $\Delta_d$, cf. \cite{cgrtv}, and the corresponding operator family is the discrete heat semigroup.

Let $s(x)=x$ and $I=[-1,1]$, say. If $\mu$ has finite moments, the standard orthonormal polynomials $\{p_n\}$ possess a three-term recurrence relation,
$$xp_n=a_{n+1}p_{n+1}+b_np_n+a_np_{n-1}.$$
In the the Schauder basis $\{\phi_n\}_n=\{p_n\}_n$ the operator can be described as a (three-diagonal) Jacobi matrix $J$.
If the Radon-Nikodym derivative of $\mu$ is positive on $(-1,1)$, by \cite[Theorem 4.5.7]{n} (see also \cite{r}) then the recurrence coefficients fulfil the asymptotics
\begin{equation}\label{sa} \lim_{n\to\infty}a_n=\frac{1}{2}\ws \ws \mbox{and} \ws \ws \lim_{n\to\infty}b_n=0,\end{equation}
and so $-2A=2(J-I)$ can be decomposed to the sum of a symmetric and a compact operator, where the symmetric part is just $\Delta_d$; in other words the rows of $-2A$ tends to the rows of $\Delta_d$. In these cases $\{W_t\}_{t\ge 0}$ called a discrete heat semigroup again, see \cite{bcf} and \cite{acl} in ultraspherical and Jacobi cases, respectively.

 Of course, by recursion for any polynomial $s(x)$ one can derive a similar recurrence relation with length $2\deg s +1$ and with $M_s=s(J)$. Here $s_M$ is the maximal element of the spectrum of the operator, that is $\max s([-1,1])$ cf. \cite[Ch. X sec. 4]{st} and on the right-hand side there is a more complicated "difference" operator. (It is a standard diffusion operator, if $1-s$ can be expressed as a polynomial of $1-x$.)

Similarly to the previous one, a recurrence relation with a $d$-diagonal matrix generates a difference operator on the right-hand side (possibly more complicated than $\Delta_d$), and by the previous computations the corresponding diffusion semigroup generates solution to the initial-value problem \eqref{h}.

In the following sections we investigate the maximal operator of the diffusion semigroup associated to multiplication operators with $d$-diagonal matrices, where $d>3$.
To derive the main results our starting point will be \eqref{6} and  the next result of Betancor et al. Before stating the theorem we need some definitions, notation.

\note

A weight on $\mathcal{I}$ is a strictly positive sequence, $w=\{w(n)\}_{n\in\mathcal{I}}$. The corresponding weighted $l^p$ spaces are
$$l^p(\mathcal{I},w)=\left\{ \{f(n)\}: \|f\|_{p,w}^p:=\sum_{m\in \mathcal{I}}|f(m)|^pw(m) <\infty\right\},$$
$1\le p<\infty$, and the weak weighted $l^1$-space is
$$l^{1,\infty}(\mathcal{I},w)=\left\{ \{f(n)\}: \|f\|_{1,\infty,w}:=\sup_{t>0}\sum_{m\in \mathcal{I}:|f(m)|>t}w(m) <\infty\right\}.$$
A weight $w=\{w(n)\}_{n\in\mathcal{I}}$ belongs to the discrete Muckenhoupt class, $A_p(\mathcal{I})$, $1\le p<\infty$ if
$$\sup_{m\le n}\frac{1}{(n-m+1)^p}\left(\sum_{k=m}^nw(k)\right)\left(\sum_{k=m}^nw(k)^{-\frac{1}{p-1}}\right)^{p-1}<\infty,$$
and it belongs to the discrete $A_1(\mathcal{I})$ class if
$$\sup_{m\le n}\frac{1}{n-m+1)}\left(\sum_{k=m}^nw(k)\right)\max w(k)^{-1}<\infty.$$

\medskip

Let $\mathbb{B}_1$ and $\mathbb{B}_2$ be Banach spaces and $\mathcal{L}(\mathbb{B}_1,\mathbb{B}_2)$ the space of bounded linear operators from $\mathbb{B}_1$ to $\mathbb{B}_2$. Let $K: (\mathbb{N}\times \mathbb{N})\setminus D \longrightarrow \mathcal{L}(\mathbb{B}_1,\mathbb{B}_2)$, where the diagonal $D$ is measurable.

\begin{defi} We say that $K$ is a local $\mathcal{L}(\mathbb{B}_1,\mathbb{B}_2)$-standard kernel, if the following conditions hold:\\
\rm {(1)}$\|K(n,m)\|_{\mathcal{L}(\mathbb{B}_1,\mathbb{B}_2)}\le \frac{C}{|n-m|},$\\
\rm {(2)}$\|K(n,m)-K(l,m)\|_{\mathcal{L}(\mathbb{B}_1,\mathbb{B}_2)}\le C\frac{|n-l|}{|n-m|^2}, \ws\ws |n-m|>2|n-l|, \ws \frac{2}{3}m\le n,l \le \frac{3}{2}m,$\\
\rm{(3)}$\|K(m,n)-K(m,l)\|_{\mathcal{L}(\mathbb{B}_1,\mathbb{B}_2)}\le C\frac{|n-l|}{|n-m|^2}, \ws\ws |n-m|>2|n-l|, \ws \frac{2}{3}m\le n,l \le \frac{3}{2}m.$\end{defi}

{\bf Theorem A.}\cite[Theorem 2.1]{bcf} {\it Let $\mathbb{B}_1$ and $\mathbb{B}_2$ be Banach spaces. Suppose that $T$ is a linear and bounded operator from $l^r_{\mathbb{B}_1}(\mathbb{N})$ into $l^r_{\mathbb{B}_2}(\mathbb{N})$ for some $1< r < \infty$, and such that there exists a local $\mathcal{L}(\mathbb{B}_1,\mathbb{B}_2)$-standard kernel $K$ such that for every finite sequence $f\in \mathbb{B}_1$
$$T(f)(n)=\sum_{m\in\mathbb{N}}K(n,m)f(m),$$
for every $n\in\mathbb{N}$, $f(n)=0$. Then,

\rm{(A1)} for every $1<p<\infty$ and $w\in A_p(\mathbb{N})$ the operator $T$ can be extended from $l^r_{\mathbb{B}_1}(\mathbb{N})\cap l^p_{\mathbb{B}_1}(\mathbb{N},w)$ to $l^p_{\mathbb{B}_1}(\mathbb{N},w)$ as a bounded operator from $l^p_{\mathbb{B}_1}(\mathbb{N},w)$ to $l^p_{\mathbb{B}_2}(\mathbb{N},w)$.

\rm{(A2)} for every $w\in A_1(\mathbb{N})$ the operator $T$ can be extended from $l^r_{\mathbb{B}_1}(\mathbb{N})\cap l^1_{\mathbb{B}_1}(\mathbb{N},w)$ to $l^1_{\mathbb{B}_1}(\mathbb{N},w)$ as a bounded operator from $l^1_{\mathbb{B}_1}(\mathbb{N},w)$ to $l^{1,\infty}_{\mathbb{B}_2}(\mathbb{N},w)$.}

\section{Discrete diffusion semigroup associated with exceptional Jacobi polynomials}

\subsection{Exceptional Jacobi polynomials} 
Introduction of exceptional orthogonal polynomials is  motivated by problems in quantum mechanics. In spite of this topic being fairly new (one of the earliest papers is \cite{gukm0}), it has a rather extended literature, see eg. \cite{ggm} and the references therein. We use the Bochner-type characterization of exceptional polynomials given in \cite{ggm}.

Classical orthogonal polynomials $\left\{P_n^{[0]}\right\}_{n=0}^\infty$ are eigenfunctions of the second order linear differential operator {with polynomial coefficients}
$$T[y]=py''+qy'+ry,$$
and its eigenvalues are denoted by $-\lambda_n$.
$T$ can be decomposed as
\begin{equation}\label{si}T=BA+\tilde{\lambda}, \ws \mbox{with} \ws A[y]=b(y'-wy), \ws B[y]=\hat{b}(y'-\hat{w}y),\end{equation}
where $b$, $w$ are rational functions and
\begin{equation}\label{s0}\hat{b}=\frac{p}{b}, \ws \ws \hat{w}=-w-\frac{q}{p}+\frac{b'}{b}.\end{equation}
Then the exceptional polynomials are the eigenfunctions of $\hat{T}$, that is the partner operator of $T$, which is
\begin{equation}\label{ka}\hat{T}[y]=(AB+\tilde{\lambda})[y]=py''+\hat{q}y'+\hat{r}y,\end{equation}
where
\begin{equation}\label{hat}\hat{q}=q+p'-2\frac{b'}{b}p, \ws \hat{r}=r+q'+wp'-\frac{b'}{b}(q+p')+\left(2\left(\frac{b'}{b}\right)^2-\frac{b''}{b}+2w'\right)p,\end{equation}
and $w$ fulfils the Riccati equation
\begin{equation}\label{w}p(w'+w^2)+qw+r=\tilde{\lambda},\end{equation}
cf. \cite[Propositions 3.5 and 3.6]{ggm}.
\eqref{si} and \eqref{ka} imply that
\begin{equation}\label{sk}\hat{T}AP_n^{[0]}=\lambda_nAP_n^{[0]},\end{equation}
so exceptional polynomials can be obtained from the classical ones by application of (finite) appropriate first order differential operator(s) to the classical polynomials. Subsequently we investigate exceptional Jacobi polynomials obtained by one Darboux transformation:
\begin{equation}\label{A} AP_n^{[0]}=b\left(P_n^{[0]}\right)'-bwP_n^{[0]}=:P_n^{[1]}.\end{equation}
The degree of $P_n^{[1]}$ is usually greater than $n$. Actually finite many ones are missing from the sequence of degrees, that is exceptional family of polynomials has finite codimension in the space of polynomials. Despite these facts, if the set of the gaps is admissible, $\left\{P_n^{[1]}\right\}_{n=0}^\infty$ is a complete orthogonal system on $I$ with respect to the weight
\begin{equation}\label{ew}W:=\frac{pw_0}{b^2},\end{equation}
where {$w_0$} is one of the classical weights, see \cite{d}, \cite{ggrm} and the references therein. To get a polynomial system, $b$ and $bw$ have to be polynomials, and in order to the moments of $W$ be finite, $b\neq 0$ on $(-1,1)$. We assume that $b>0$ on $(-1,1)$.

After this general summary we introduce the exceptional Jacobi polynomials. Our starting point is the classical Jacobi system. We mostly follow the notation of \cite{sz}.
$$w^{\alpha,\beta}=(1-x)^\alpha(1+x)^\beta$$
$$p_k^{\alpha,\beta}=\frac{P_k^{\alpha,\beta}}{\varrho_k^{\alpha,\beta}},$$
where
\begin{equation}\label{ro}\left(\varrho_k^{\alpha,\beta}\right)^2= \frac{2^{\varrho}\Gamma(k+\alpha+1)\Gamma(k+\beta+1)}{(2k+\varrho)\Gamma(k+1)\Gamma(k+\varrho)}.\end{equation}
with
$$\varrho:=\alpha+\beta+1.$$
\begin{equation}\tilde{P}_n=\tilde{P}_n^{\alpha,\beta}=\frac{1}{\sigma_n} P_n^{[1]},\end{equation}
where
$$\sigma_n:=\sigma_n^{\alpha,\beta}=\|P_n^{[1]}\|_{W,2},$$
and
$$P_n^{[1]}=P_n^{\alpha, \beta, [1]}=b(p_n^{\alpha,\beta})'-bwp_n^{\alpha,\beta}.$$
Subsequently we assume that the admissibility condition mentioned above fulfils, that is the system is complete.

Because finitely many degrees are missing from the sequence of degrees, exceptional orthogonal polynomials do not fulfil three-term recurrence formulae, and it can happen that $x\tilde{P}_n$ can be expressed as an infinite series with respect to $\{\tilde{P}_n\}$. Fortunately it is proved that similarly to the standard cases, there are finite recurrence relations with certain polynomials, more precisely with the notation $\frac{p}{b}=\frac{\tilde{p}}{\tilde{b}}$ (see \eqref{ew}) if $\tilde{b}$  is a divisor of $s'$, then $M_s$ is a $2\mathrm{deg}s+1$-diagonal infinite matrix in the basis of $\{\tilde{P}_n\}$.

If $s_M-s$ has a simple zero at $1$, the situation is rather similar to the one discussed in \cite{acl}. Below we assume that
\begin{equation}\label{b1} b(1)=0 \ws \mbox{ and} \ws b'(1)\neq 0 ,\end{equation}
and
define $s=:Q$ as a primitive function of $b$.
\begin{equation}\label{q}Q(x):= \int^x b.\end{equation}
Of course, the constant term of $Q$ can be chosen. With $s(x):=Q(x)$ we get the almost simplest recurrence relation
\begin{equation}\label{R}Q\tilde{P}_n=\sum_{k=-L}^Lu_{n,k}\tilde{P}_{n+k},\end{equation}
where $L$ is the degree of $Q$, see \cite{o}. 
Multiplication operator with respect to exceptional Jacobi polynomials is examined in \cite{h} and \cite{h2}. The matrix of the corresponding multiplication operator in $\{\tilde{P}_n\}$ basis is
\begin{equation}\label{mQ}M_{Q}=\left[\begin{array}{cccccccc}u_{0,0}&u_{0,1}&\dots&\dots&u_{0,L}&0&0&\dots\\u_{1,-1}&u_{1,0}&\dots&\dots&u_{1,L-1}&u_{1,L}&0&\dots\\ \vdots&\vdots&\dots&\ddots&\vdots&\dots&\vdots&\dots\\u_{L,-L}&u_{L,-L+1}&\dots&u_{L,0}&\dots&\dots&u_{L,L}&0\\ 0&u_{L+1,-L}&\dots&\vdots&\dots&\vdots&\dots&u_{L+1,L}\\ \vdots&0&\dots&u_{L+j,-L}& \dots&\vdots&\dots&\dots\end{array}\right].\end{equation}
It can be easily seen that $M_Q$ is symmetric since
\begin{equation}\label{ak}u_{k,j}=\int_{-1}^1Q\tilde{P}_k\hat{P}_{k+j}W^2=\int_{-1}^1Q\tilde{P}_{k+j}\tilde{P}_{(k+j)-j}W^2=u_{k+j,-j}.\end{equation}
Furthermore the coefficients in \eqref{R} fulfil the symmetric limit relation
\begin{equation}\label{UU}\lim_{n\to\infty}u_{n,j}=:U_{|j|},\end{equation}
where $U_{|j|}$ depends on the polynomial $b$, see \cite[Proposition 3.4]{h}.

\subsection{The discrete diffusion semigroup}
With this $Q$ the next (symmetric) kernel can be defined.
\begin{equation}K^{\alpha,\beta,e}_t(n,m):=\int_{-1}^1e^{-(Q(1)-Q(x))t}\tilde{P}_n(x)\tilde{P}_m(x)W(x)dx,\end{equation}
and the operator
\begin{equation}W^{\alpha,\beta,e}_tf(n):=\sum_{m=0}^\infty f(m)K^{\alpha,\beta,e}_t(n,m).\end{equation}

\medskip

According to the results of the previous section $\{W^{\alpha,\beta,e}_t\}_{t\ge 0}$ is the discrete diffusion semigroup associated with exceptional Jacobi polynomials and $W_tf(n)$ is a solution to \eqref{h}.\\
 Unlike the ultraspherical and the Jacobi cases, the operators are not positivity preserving.\\ 
Before stating the main result of this section, let us have an example. Let
$$b(x):=(1-x)P_1^{(-\frac{3}{2}, \frac{1}{2})}=\frac{x^2}{2}-\frac{3}{2}x+1.$$
In view of \cite[(A),(B), (89)]{gumm} it is an appropriate choice. By \cite[(58)-(60)]{gumm}
$$w(x)=\frac{1}{2}\frac{P_1^{(-\frac{5}{2}, -\frac{1}{2})}(x)}{(1-x)P_1^{(-\frac{3}{2}, \frac{1}{2})}(x)},$$
and
$$P_n^{[1]}(x)=(1-x)P_1^{(-\frac{3}{2}, \frac{1}{2})}(x)(p_n^{(\frac{3}{2}, \frac{1}{2})}(x))'-\frac{1}{2}P_1^{(-\frac{5}{2}, -\frac{1}{2})}(x)p_n^{(\frac{3}{2}, \frac{1}{2})}(x).$$
Let
$$Q(x)=\frac{x^3}{6}-\frac{3}{4}x^2+x.$$
Then, by \cite[(3.20)]{h}
$$U_0=-\frac{3}{8}, \ws \ws U_1=\frac{9}{16}, \ws \ws U_2=-\frac{3}{16}, \ws \ws U_3=\frac{1}{48}.$$
Thus $48A=48(M_Q-Q(1)I)$ can be decomposed to the sum of a symmetric and a compact operator. The rows of the symmetric one are
$$\dots, 0,1,-9,27,-38,27,-9,1, 0,\dots.$$

Thus in this case $\{W_t^{\alpha,\beta,e}\}_{t\ge 0}$ is a diffusion semigroup associated to the initial-value problem
\begin{equation}\left\{\begin{array}{ll}\frac{\partial u(n,t)}{\partial t}=L_du(n,t),\\ u(n,0)=f(n)\end{array}\right.,\end{equation}
cf.\eqref{h}. $L_d=L_{d,s}+L_{d,c}$, where $L_{d,s}$ is the discrete version of
$$L[u]=\frac{1}{48}\left(\frac{\partial^6 u}{\partial x^6}-3\frac{\partial^4 u}{\partial x^4}\right),$$
because $\Delta x=1$.

\subsection{The maximal operator}
Now we are in position to state the main result of this section about the maximal operator \eqref{max}.

\begin{theorem}\label{t1} Supposing the assumptions \eqref{b1} and \eqref{q} are satisfied if $\alpha>\frac{3}{2}$, $\beta\ge-\frac{1}{2}$, then the maximal operator of the discrete diffusion semigroup associated with exceptional Jacobi polynomials fulfils that\\
\rm{(1)} if $1<p<\infty$ and $w\in A_p(\mathbb{N})$, then for all $f\in l^2(\mathbb{N})\cap l^p(\mathbb{N},w)$
$$\|W^{\alpha,\beta,e}_*f\|_{p,w}\le C\|f\|_{p,w},$$
where $C$ is a constant independent of $f$. That is the operator $W^{\alpha,\beta,e}_*$ can be extended uniquely to a bounded operator from $l^p(\mathbb{N},w)$ into itself.\\
\rm{(2)} If $w\in A_p(\mathbb{N})$, then for all $f\in l^2(\mathbb{N})\cap l^1(\mathbb{N},w)$
$$\|W^{\alpha,\beta,e}_*f\|_{(1,\infty),w}\le C\|f\|_{1,w},$$
where $C$ is a constant independent of $f$. That is the operator $W^{\alpha,\beta,e}_*$ can be extended uniquely to a bounded operator from $l^1(\mathbb{N},w)$ into $l^{1,\infty}(\mathbb{N},w)$.
\end{theorem}

\medskip

\remark
The unpleasant requirement, $\alpha>\frac{3}{2}$, is necessary because $Q(1)-Q(x)$ has a double zero at $x=1$.

\medskip

According to \cite[under (38)]{bcf} or \cite[(21)]{acl}, in order to prove Theorem \ref{t1} it is enough to prove the next lemma.

\medskip

\begin{lemma}\label{le1}
\begin{equation}\label{l2}|K^{\alpha,\beta,e}_t(n,n)|\le C,\end{equation}
Let $\alpha >\frac{3}{2}$, $\beta \ge-\frac{1}{2}$. Then
\begin{equation}\label{l1}|K^{\alpha,\beta,e}_t(n,m)| \le \frac{C}{|n-m|}, \ws \ws \mbox{if} \ws \ws n\neq m ,\end{equation}
and
\begin{equation}\label{l3}|K^{\alpha,\beta,e}_t(n+1,m)-K^{\alpha,\beta,e}_t(n,m)| \le \frac{C}{|n-m|^2}, \ws \ws \mbox{if} \ws \ws n\neq m,m\pm 1,\ws \frac{m}{2}<n<2m, \end{equation}
where $C$ is a constant which may be different at each occurrence, and is independent of $t$, $n$ and $m$.
\end{lemma}

\medskip

Before the proof we introduce some notations for sake of convenience. First let
$$\overline{Q}(x):=(Q(1)-Q(x)).$$
Below we use
\begin{equation}I_{n,m}^{\alpha,\beta,\gamma,\delta,\mu,\nu}[f,t]:=\int_{-1}^1e^{-\overline{Q}t}p_n^{\alpha,\beta}p_m^{\gamma,\delta}f_tw^{\mu,\nu}dx,\end{equation}
\begin{equation}D_{n,m}^{\alpha,\beta,\gamma,\delta,\mu,\nu}[f,t]:=\int_{-1}^1e^{-\overline{Q}t}(p_n^{\alpha,\beta}-p_{n-1}^{\alpha,\beta})p_m^{\gamma,\delta}f_tw^{\mu,\nu}dx,\end{equation}
for $\gamma=\mu=\alpha$, $\delta=\nu=\beta$
\begin{equation}\label{Ihat}\hat{I}_{n,m}^{\alpha,\beta}[f,t]:=\int_{-1}^1e^{-\overline{Q}t}p_n^{\alpha,\beta}p_m^{\alpha,\beta}f_t(x)w^{\alpha,\beta}dx,\end{equation}
and similarly
\begin{equation}\label{Dhat}\hat{D}_{n,m}^{\alpha,\beta}[f,t]=D_{n,m}^{\alpha,\beta,\alpha,\beta,\alpha,\beta}[f,t],\end{equation}
the operator
\begin{equation}Lf:=tbf_t+f'_t,\end{equation}
and the constant
$$\sigma=\gamma+\delta+1.$$
We use the notation $f'_t:= \frac{\partial}{\partial x}f_t(x)$.

\medskip

\proof

According to the construction we can express the kernel function of the operator by classical Jacobi polynomials.
$$K^{\alpha,\beta,e}_t(n,m)=\frac{1}{\sigma_n \sigma_m}\int_{-1}^1e^{-\overline{Q}(x)t}$$ $$\times\left(b(x)(p_n^{\alpha,\beta})'(x)-(bw)(x)p_n^{\alpha,\beta}(x)\right)\left(b(x)(p_m^{\alpha,\beta})'(x)-(bw)(x)p_m^{\alpha,\beta}(x)\right)\frac{p(x)w^{\alpha,\beta}(x)}{b^2(x)}dx$$
$$=\frac{1}{\sigma_n \sigma_m}\int_{-1}^1e^{-\overline{Q}t}\left((p_n^{\alpha,\beta})'(p_m^{\alpha,\beta})'w^{\alpha+1,\beta+1}+pw^2p_n^{\alpha,\beta}p_m^{\alpha,\beta}w^{\alpha,\beta}-pw\left(p_n^{\alpha,\beta}p_m^{\alpha,\beta}\right)'w^{\alpha,\beta}\right)dx.$$
Integrating by parts in the third term and considering that
\begin{equation}\label{d}(p_n^{\alpha,\beta})'=\sqrt{n(n+\varrho)}p_{n-1}^{\alpha+1,\beta+1},\end{equation}
(see \cite[(4.21.7)]{sz})
$$K^{\alpha,\beta,e}_t(n,m)$$ $$=\frac{1}{\sigma_n \sigma_m}\int_{-1}^1e^{-\overline{Q}t}\left((p_n^{\alpha,\beta})'(p_m^{\alpha,\beta})'w^{\alpha+1,\beta+1}+p_n^{\alpha,\beta}p_m^{\alpha,\beta}(btpw+p(w^2+w')+qw)w^{\alpha,\beta}\right)dx.$$

By \eqref{Ihat} and considering \eqref{w} (with $r\equiv 0$), we have
\begin{equation}\label{K}K^{\alpha,\beta,e}_t(n,m)=:K_1(n,m)+K_2(n,m)\end{equation}
$$=\frac{\sqrt{n(n+\varrho)}\sqrt{m(m+\varrho)}}{\sigma_n \sigma_m}\hat{I}_{n-1,m-1}^{\alpha+1,\beta+1}[1,t]+\frac{1}{\sigma_n \sigma_m}\hat{I}_{n,m}^{\alpha,\beta}[tbpw+\tilde{\lambda},t].$$
Now we turn to the proof of the statements of the lemma.

Recalling that $b$ is bounded on $(-1,1)$, \eqref{l2} follows from orthogonality.\\
By the symmetry of the kernel we can choose $n>m$, say. Following the chain of ideas of \cite{acl} first we apply that
\begin{equation}\label{dd}(p_n^{\alpha,\beta}w^{\alpha,\beta})'=-\sqrt{(n+1)(n+\varrho-1)}p_{n+1}^{\alpha-1,\beta-1}w^{\alpha-1,\beta-1},\end{equation}
(see \cite[(4.10.1)]{sz}) with $\alpha,\beta,\gamma,\delta,\mu,\nu >-1$, and then \eqref{d}. Thus
$$I_{n,m}^{\alpha,\beta,\gamma,\delta,\mu,\nu}[f,t]=-\frac{1}{\sqrt{n(n+\varrho)}}\int_{-1}^1e^{-\overline{Q}t}(p_{n-1}^{\alpha+1,\beta+1}w^{\alpha+1,\beta+1})'p_m^{\gamma,\delta}f_tw^{\mu-\alpha,\nu-\beta}dx$$
$$=\frac{1}{\sqrt{n(n+\varrho)}}\int_{-1}^1e^{-\overline{Q}t}p_{n-1}^{\alpha+1,\beta+1}w^{\alpha+1,\beta+1}\left(\sqrt{m(m+\sigma)}p_{m-1}^{\gamma+1,\delta+1}f_tw^{\mu-\alpha,\nu-\beta}\right.$$ $$\left.+p_m^{\gamma,\delta}f'_tw^{\mu-\alpha,\nu-\beta}+p_m^{\gamma,\delta}f_tw^{\mu-\alpha-1,\nu-\beta-1}(\nu-\beta-\mu+\alpha-(\nu-\alpha-\beta+\mu)x)\right.$$ $$\left.+btp_m^{\gamma,\delta}f_tw^{\mu-\alpha,\nu-\beta}\right)dx=(*).$$
Proceeding in the same way with the first term we have
$$(*)=\frac{m(m+\sigma)}{n(n+\varrho)} I_{n,m}^{\alpha,\beta,\gamma,\delta,\mu,\nu}[f,t]+\frac{1}{\sqrt{n(n+\varrho)}}\int_{-1}^1e^{-\overline{Q}t}p_{n-1}^{\alpha+1,\beta+1}p_m^{\gamma,\delta}\left(f'_tw^{\mu+1,\nu+1}\right.$$ $$\left.+f_tw^{\mu,\nu}(\nu-\beta-\mu+\alpha-(\nu-\alpha-\beta+\mu)x)+btf_tw^{\mu+1,\nu+1}\right)dx $$ $$- \frac{\sqrt{m(m+\sigma)}}{n(n+\varrho)}\int_{-1}^1e^{-\overline{Q}t}p_n^{\alpha,\beta}p_{m-1}^{\gamma+1,\delta+1}\left(w^{\mu,\nu}(\nu-\delta-\mu+\gamma-(\nu-\delta-\gamma+\mu)x)f_t\right.$$ $$\left.+f'_tw^{\mu+1,\nu+1}+f_tw^{\mu+1,\nu+1}bt\right)dx.$$
So
\begin{equation}\label{I}I_{n,m}^{\alpha,\beta,\gamma,\delta,\mu,\nu}[f,t]\end{equation} $$=\frac{1}{n-m}\left(\frac{\sqrt{n(n+\varrho)}}{n+m+\varrho+\frac{m}{n-m}(\varrho-\sigma)}I_{n-1,m}^{\alpha+1,\beta+1,\gamma,\delta,\mu+1,\nu+1}[Lf,t]\right.$$ $$\left.- \frac{\sqrt{m(m+\sigma)}}{n+m+\varrho+\frac{m}{n-m}(\varrho-\sigma)}I_{n,m-1}^{\alpha,\beta,\gamma+1,\delta+1,\mu+1,\nu+1}[Lf,t]\right.$$ $$\left.+\frac{\sqrt{n(n+\varrho)}}{n+m+\varrho+\frac{m}{n-m}(\varrho-\sigma)}I_{n-1,m}^{\alpha+1,\beta+1,\gamma,\delta,\mu,\nu}[(\nu-\beta-\mu+\alpha-(\nu-\alpha-\beta+\mu)x)f_t,t]\right.$$ $$\left.- \frac{\sqrt{m(m+\sigma)}}{n+m+\varrho+\frac{m}{n-m}(\varrho-\sigma)}I_{n,m-1}^{\alpha,\beta,\gamma+1,\delta+1,\mu,\nu}[(\nu-\delta-\mu+\gamma-(\nu-\delta-\gamma+\mu)x)f_t,t]\right).$$
To prove \eqref{l1} considering \eqref{I} first we deal with the first term of \eqref{K}.
\begin{equation}\label{nm}\hat{I}_{n-1,m-1}^{\alpha+1,\beta+1}[1,t]=\frac{1}{n-m}\left(\frac{\sqrt{(n-1)(n+\varrho+1)}}{n+m+\varrho}I_{n-2,m-1}^{\alpha+2,\beta+2,\alpha+1,\beta+1,\alpha+2,\beta+2}[tb,t]\right.$$ $$\left.- \frac{\sqrt{(m-1)(m+\varrho+1)}}{n+m+\varrho}I_{n-1,m-2}^{\alpha+1,\beta+1,\alpha+2,\beta+2,\alpha+2,\beta+2}[tb,t]\right).\end{equation}
Thus we have to estimate the two integrals $I_{n-2,m-1}^{\alpha+2,\beta+2,\alpha+1,\beta+1,\alpha+2,\beta+2}[tb,t]$ and $I_{n-1,m-2}^{\alpha+1,\beta+1,\alpha+2,\beta+2,\alpha+2,\beta+2}[tb,t]$. Since the computations are the same we deal with the second one, say.
We apply the next norm estimation
\begin{equation}\label{n}\left\|p_{n}^{\alpha,\beta}w^{\frac{\alpha}{2}+\frac{1}{4}, \frac{\beta}{2}+\frac{1}{4}}\right\|_\infty \le C,\end{equation}
where $C$ is an absolute constant (see \cite[(8.21.10)]{sz}). Recalling again that $b>0$ on $(-1, 1)$ we have
$$I_{n-1,m-2}^{\alpha+1,\beta+1,\alpha+2,\beta+2,\alpha+2,\beta+2}[tb,t] $$ $$\le \|p_{n-1}^{\alpha+1,\beta+1}w^{\frac{\alpha}{2}+\frac{3}{4},\frac{\beta}{2}+\frac{3}{4}}\|_\infty\|p_{m-2}^{\alpha+2,\beta+2}w^{\frac{\alpha}{2}+\frac{5}{4},\frac{\beta}{2}+\frac{5}{4}}\|_\infty\int_{-1}^1e^{-\overline{Q}t}btdx<C,$$
 because the integral is uniformly bounded in $t$. The coefficients of the two integrals in \eqref{nm} are less than 1, thus
\begin{equation}\label{nm1}\left|\hat{I}_{n-1,m-1}^{\alpha+1,\beta+1}[1,t]\right|<\frac{C}{|n-m|},\end{equation}
and so
\begin{equation}\label{K1}|K_1(n,m)|<\frac{C}{|n-m|}.\end{equation}
In view of \eqref{I}
\begin{equation}\label{nm2}K_2(n,m)=\frac{1}{\sigma_n\sigma_m}\hat{I}_{n,m}^{\alpha,\beta}[tbpw+\tilde{\lambda},t]=\frac{1}{\sigma_n\sigma_m(n-m)}\end{equation} $$\times\left(\frac{\sqrt{n(n+\varrho)}}{n+m+\varrho}I_{n-1,m}^{\alpha+1,\beta+1,\alpha,\beta,\alpha+1,\beta+1}[t^2b^2pw+t(b\tilde{\lambda}+(bpw)'),t]\right.$$ $$\left.- \frac{\sqrt{m(m+\varrho)}}{n+m+\varrho}I_{n,m-1}^{\alpha,\beta,\alpha+1,\beta+1,\alpha+1,\beta+1}[t^2b^2pw+t(b\tilde{\lambda}+(bpw)'),t]\right).$$
The two integrals can be handled on the same way again. We decompose the first one, say, to two terms; $I_1$ contains the members of the first order in $t$, and $I_2$ the ones of the second order.
$$I_1:=I_{n-1,m}^{\alpha+1,\beta+1,\alpha,\beta,\alpha+1,\beta+1}[t(b\tilde{\lambda}+(bpw)'),t]$$ $$=\int_{-1}^1e^{-\overline{Q}t}\overline{Q}t\frac{(1-x)^2}{\overline{Q}}(b\tilde{\lambda}+(bpw)')p_{n-1}^{\alpha+1,\beta+1}p_m^{\alpha,\beta}w^{\alpha-1,\beta+1}dx.$$
Let us recall that $bpw$ is a polynomial, $\overline{Q}$ has a double zero at $x=1$ and as $b>0$ on $(-1,1)$, $Q$ is monotone there. Considering that $e^{-u}u$ is uniformly bounded, we have
$$|I_1|\le C\|p_{n-1}^{\alpha+1,\beta+1}w^{\frac{\alpha}{2}+\frac{3}{4},\frac{\beta}{2}+\frac{3}{4}}\|_\infty\int_{-1}^1|p_m^{\alpha,\beta}|w^{\frac{\alpha}{2}-\frac{7}{4},\frac{\beta}{2}+\frac{1}{4}}dx\le Cm^2,$$
where the last integral is convergent by the assumption on $\alpha$ and $\beta$ and the estimation follows from
\begin{equation}\label{ib}\int_0^1|(1-x)^\mu p_n^{\alpha,\beta}(x)|dx\sim n^{\alpha+\frac{1}{2}-2\mu-2},\end{equation}
if $2\mu<\alpha-\frac{3}{2}$, see \cite[(7.34.1)]{sz}.
$$I_2:=I_{n-1,m}^{\alpha+1,\beta+1,\alpha,\beta,\alpha+1,\beta+1}[t^2b^2pw,t]$$ $$=\int_{-1}^1e^{-\overline{Q}t}\overline{Q}^2t^2\frac{b^2(1-x)^2}{\overline{Q}^2}pwp_{n-1}^{\alpha+1,\beta+1}p_m^{\alpha,\beta}w^{\alpha-1,\beta+1}dx.$$
This estimation can be finished as the previous one: $pw$ and the fraction and $e^{-u}u^2$ are bounded and then we get the same integral to estimate. Thus
$$|I_1|+|I_2|<C m^2,$$ and taking into consideration that
\begin{equation}\label{s}\sigma_k=\sqrt{k(k+\varrho)+\tilde{\lambda}},\end{equation}
see\cite[(42)]{h2}
\begin{equation}\label{nm2}\frac{1}{\sigma_n \sigma_m}|\hat{I}_{n,m}^{\alpha,\beta}[tbpw+\tilde{\lambda},t]|\le C \frac{m^2}{mn|n-m|}\le \frac{C}{|n-m|}.\end{equation}
Finally considering \eqref{K}, \eqref{s}, \eqref{nm1} and \eqref{nm2}, \eqref{l1} is proved.

To prove \eqref{l3} we iterate \eqref{I} once more, that is if $n-1>m$, say
\begin{equation}\label{I2}\hat{I}_{n,m}^{\alpha,\beta}[f,t]\end{equation} $$=\frac{\sqrt{n(n+\varrho)(n-1)(n+\varrho+1)}I_{n-2,m}^{\alpha+2,\beta+2, \alpha,\beta,\alpha+2,\beta+2}[L^2f,t]}{(n-m)(n-m-1)(n+m+\varrho)(n+m+\varrho+1+\frac{2m}{n-m-1})}$$
$$+\frac{\sqrt{m(m+\varrho)(m-1)(m+\varrho+1)}I_{n,m-2}^{\alpha,\beta, \alpha+2,\beta+2,\alpha+2,\beta+2}[L^2f,t]}{(n-m)(n-m+1)(n+m+\varrho)(n+m+\varrho-1-\frac{2(m-1)}{n-m+1})}$$
$$-\frac{\sqrt{n(n+\varrho)m(m+\varrho)}I_{n-1,m-1}^{\alpha+1,\beta+1, \alpha+1,\beta+1,\alpha+2,\beta+2}[L^2f,t]}{(n-m)(n-m-1)(n+m+\varrho)(n+m+\varrho+1+\frac{2m}{n-m-1})}$$
$$-\frac{\sqrt{m(m+\varrho)n(n+\varrho)}I_{n-1,m-1}^{\alpha+1,\beta+1, \alpha+1,\beta+1,\alpha+2,\beta+2}[L^2f,t]}{(n-m)(n-m+1)(n+m+\varrho)(n+m+\varrho-1-\frac{2(m-1)}{n-m+1})}$$
$$+\frac{\sqrt{n(n+\varrho)m(m+\varrho)}\hat{I}_{n-1,m-1}^{\alpha+1,\beta+1}[2xLf,t]}{(n-m)(n-m-1)(n+m+\varrho)(n+m+\varrho+1+\frac{2m}{n-m-1})}$$
$$+\frac{\sqrt{n(n+\varrho)m(m+\varrho)}\hat{I}_{n-1,m-1}^{\alpha+1,\beta+1}[2xLf,t]}{(n-m)(n-m+1)(n+m+\varrho)(n+m+\varrho-1-\frac{2(m-1)}{n-m+1})}.$$
(It is clear that the estimation of the first two terms, the two terms in the middle and the last two terms are the same. )

In view of \eqref{K} we start with the estimation
$$|K_1(n+1,m)-K_1(n,m)|\le O\left(\frac{1}{n^3}\right)\left|\hat{I}_{n,m-1}^{\alpha+1,\beta+1}[1,t]\right| +O(1)\left|\hat{D}_{n,m-1}^{\alpha+1,\beta+1}[1,t]\right|.$$
As $\hat{I}_{n,m-1}^{\alpha+1,\beta+1}[1,t]$ is uniformly bounded, the estimation in \eqref{l3} is obvious for the first term. To estimate the second one we make the next observation.

By the assumption on $n,m$ the denominators of the coefficients of the integrals are of $O\left(\frac{1}{(n-m)^2n^2}\right)$ thus according to \eqref{I2}
\begin{equation}\label{D}\left|\hat{D}_{n,m}^{\alpha,\beta}[L^2f,t]\right|\end{equation}
$$\le \frac{C}{|n-m|^2}\left(\left|D_{n-1,m}^{\alpha+2,\beta+2, \alpha,\beta,\alpha+2,\beta+2}[L^2f,t]\right|+\left|D_{n+1,m-2}^{\alpha,\beta, \alpha+2,\beta+2,\alpha+2,\beta+2}[L^2f,t]\right|\right.$$ $$\left. +\left|D_{n,m-1}^{\alpha+1,\beta+1, \alpha+1,\beta+1,\alpha+2,\beta+2}[L^2f,t]\right|+\left|\hat{D}_{n,m-1}^{\alpha+1,\beta+1}[2xLf,t]\right|\right)$$ $$+\frac{C}{|n-m|^2n}\left(\left|I_{n-1,m}^{\alpha+2,\beta+2, \alpha,\beta,\alpha+2,\beta+2}[L^2f,t]\right|+\left|I_{n+1,m-2}^{\alpha,\beta, \alpha+2,\beta+2,\alpha+2,\beta+2}[L^2f,t]\right|\right.$$ $$\left.+\left|I_{n,m-1}^{\alpha+1,\beta+1, \alpha+1,\beta+1,\alpha+2,\beta+2}[L^2f,t]\right|+\left|\hat{I}_{n,m-1}^{\alpha+1,\beta+1}[2xLf,t]\right|\right)$$ $$\frac{C}{|n-m|^2}\mathcal{D}+\mathcal{R}.$$
To estimate the $"D"$ terms we use the following formula
\begin{equation}\label{dif}p_n^{\alpha+1,\beta+1}-p_{n-1}^{\alpha+1,\beta+1}=\frac{\alpha+1}{n-1}p_{n-1}^{\alpha+1,\beta+1}\end{equation} $$+\left(1-\frac{\varrho_n^{\alpha+1,\beta+1}}{\varrho_{n-1}^{\alpha+1,\beta+1}}\right)p_n^{\alpha+1,\beta+1}-\frac{2n+\varrho+1}{2n}\frac{\varrho_{n-1}^{\alpha+2,\beta+1}}{\varrho_{n-1}^{\alpha+1,\beta+1}}(1-x)p_{n-1}^{\alpha+2,\beta+1},$$
see \cite[(4.5.4)]{sz}.
As $\left(1-\frac{\varrho_n^{\alpha+1,\beta+1}}{\varrho_{n-1}^{\alpha+1,\beta+1}}\right)=O\left(\frac{1}{n}\right)$, each $"D"$ term can be divided to $\frac{1}{n}"I"$ terms which can be added to $\mathcal{R}$, and terms wich contains an extra $(1-x)$. For instance let us start with
$$\left|D_{n,m-1}^{\alpha+1,\beta+1, \alpha+1,\beta+1,\alpha+2,\beta+2}[L^21,t]\right|$$ $$\leq  \frac{C}{n}\left|\int_{-1}^1e^{-\overline{Q}t}(b^2t^2+tb')p_{n^*}^{\alpha+1,\beta+1}p_{m}^{\alpha+1,\beta+1}w^{\alpha+2,\beta+2}dx\right|$$ $$+C\int_{-1}^1e^{-\overline{Q}t}b^2t^2\left|p_{n-1}^{\alpha+2,\beta+1}p_{m}^{\alpha+1,\beta+1}\right|w^{\alpha+3,\beta+2}dx$$ $$+C\int_{-1}^1e^{-\overline{Q}t}|b'|t\left|p_{n-1}^{\alpha+2,\beta+1}p_{m}^{\alpha+1,\beta+1}\right|w^{\alpha+3,\beta+2}dx$$ $$=J_1+J_2+J_3,$$
where $n^*=n$ or $n^*=n-1$.
$J_1$ can be added to $\mathcal{R}$, and
$$J_2\le C\|p_{n-1}^{\alpha+2,\beta+1}p_{m}^{\alpha+1,\beta+1}w^{\alpha+2,\beta+\frac{3}{2}}\|_\infty \int_{-1}^1e^{-\overline{Q}t}b^2(1-x)t^2dx$$ $$\le C\int_0^1e^{-\overline{Q}t}b^3t^2dx\le C,$$
where the last estimation comes from an integration by parts. Similarly
$$J_3\le C \int_{-1}^1e^{-\overline{Q}t}btdx\le C.$$
The first two terms in \eqref{D} can be handled on the same way - we estimate the firs one, say.
$$\left|D_{n-1,m}^{\alpha+2,\beta+2, \alpha,\beta,\alpha+2,\beta+2}[L^21,t]\right|$$ $$\le\frac{C}{n}\left|\int_{-1}^1e^{-\overline{Q}t}(b^2t^2+tb')p_{(n-1)^*}^{\alpha+2,\beta+2}p_{m}^{\alpha,\beta}w^{\alpha+2,\beta+2}dx\right|$$ $$+C\int_{-1}^1e^{-\overline{Q}t}b^2t^2\left|p_{n-2}^{\alpha+3,\beta+2}p_{m}^{\alpha,\beta}\right|w^{\alpha+3,\beta+2}dx$$ $$+C\int_{-1}^1e^{-\overline{Q}t}|b'|t\left|p_{n-2}^{\alpha+3,\beta+2}p_{m}^{\alpha,\beta}\right|w^{\alpha+3,\beta+2}dx$$ $$=J_4+J_5+J_6.$$
As above, $J_4$ can be added to $\mathcal{R}$, and
$$J_5\le C\|p_{n-2}^{\alpha+3,\beta+2}p_{m}^{\alpha,\beta}w^{\alpha+2,\beta+1}\|_\infty \int_{-1}^1e^{-\overline{Q}t}b^2(1-x)t^2dx\le C.$$
Similarly
$$J_6\le C \int_{-1}^1e^{-\overline{Q}t}btdx\le C.$$
To finish this part we estimate
$$\left|\hat{D}_{n,m-1}^{\alpha+1,\beta+1}[2xL1,t]\right|$$ $$\le \frac{C}{n}\left|\int_{-1}^1e^{-\overline{Q}t}2xbtp_{(n-1)^*}^{\alpha+1,\beta+1}p_{m-1}^{\alpha+1,\beta+1}w^{\alpha+1,\beta+1}dx\right|$$ $$+C\int_{-1}^1e^{-\overline{Q}t}2xbt\left|p_{n-1}^{\alpha+2,\beta+1}p_{m-1}^{\alpha+1,\beta+1}\right|w^{\alpha+2,\beta+1}dx=J_7+J_8.$$
Again $J_7$ can be added to $\mathcal{R}$, and
$$J_8\le C\|p_{n-1}^{\alpha+2,\beta+1}p_{m-1}^{\alpha+1,\beta+1}w^{\alpha+2,\beta+1}\|_\infty \int_{-1}^1e^{-\overline{Q}t}btdx\le C.$$
Now we turn to the estimation of terms in $\mathcal{R}$. For instance let us see the first one.
$$\frac{C}{n}\left|I_{n-1,m}^{\alpha+2,\beta+2, \alpha,\beta,\alpha+2,\beta+2}[L^21,t]\right|\le \frac{C}{n}\int_{-1}^1e^{-\overline{Q}t}b^2t^2\left|p_{n-1}^{\alpha+2,\beta+2}p_{m}^{\alpha,\beta}\right|w^{\alpha+2,\beta+2}dx$$ $$+\frac{C}{n}\int_{-1}^1e^{-\overline{Q}t}b't\left|p_{n-1}^{\alpha+2,\beta+2}p_{m}^{\alpha,\beta}\right|w^{\alpha+2,\beta+2}dx=I_3+I_4.$$
In view of \eqref{n} and \eqref{ib}
$$I_3\le \frac{C}{n}\int_{-1}^1e^{-\overline{Q}t}b^2(1-x)^2t^2\left|p_{n-1}^{\alpha+2,\beta+2}p_{m}^{\alpha,\beta}\right|w^{\alpha,\beta+2}dx$$ $$\le \frac{C}{n}\int_{-1}^1\left|p_{n-1}^{\alpha+2,\beta+2}\right|w^{\frac{\alpha}{2}-\frac{1}{4},\frac{\beta}{2}+\frac{3}{4}}dx\le \frac{C}{n}n\le C.$$
Similarly by \eqref{n}
$$I_4\le \frac{C}{n}\int_{-1}^1e^{-\overline{Q}t}(1-t)^2t\left|p_{n-1}^{\alpha+2,\beta+2}\right|w^{\frac{\alpha}{2}-\frac{1}{4},\frac{\beta}{2}+\frac{3}{4}}dx\le C.$$
The previous computations ensures that the remainder terms of $\mathcal{R}$ can be estimated similarly.

To continue the proof of \eqref{l3} in view of \eqref{s} we have
\begin{equation}\label{k2}|K_2(n+1,m)-K_2(n,m)|\end{equation} $$\le O\left(\frac{1}{mn^2}\right)\left|\hat{I}_{n+1,m}^{\alpha,\beta}[tbpw+\tilde{\lambda},t]\right| + O\left(\frac{1}{mn}\right)\left|\hat{D}_{n+1,m}^{\alpha,\beta}[tbpw+\tilde{\lambda},t]\right|.$$
$$\left|\hat{I}_{n+1,m}^{\alpha,\beta}[tbpw+\tilde{\lambda},t]\right|\le |\tilde{\lambda}|\|e^{-\overline{Q}t}\|_\infty \|p_{n+1}^{\alpha,\beta}\|_{w^{\alpha,\beta},2}\|p_{m}^{\alpha,\beta}\|_{w^{\alpha,\beta},2}$$ $$+\|pw\|_\infty \|p_{n+1}^{\alpha,\beta}p_{m}^{\alpha,\beta}w^{\alpha,\beta}\|_\infty\int_{-1}^1e^{-\overline{Q}t}btdx\le C\sqrt{mn}.$$
Recalling that $n\sim m$
$$O\left(\frac{1}{mn^2}\right)\left|\hat{I}_{n+1,m}^{\alpha,\beta}[tbpw+\tilde{\lambda},t]\right|=O\left(\frac{1}{n^2}\right).$$
To estimate the second term we proceed as above. Taking into consideration the properties of $b$ and $w$
\begin{equation}\label{L}L^2(tbpw+\tilde{\lambda})=t^3b^3pw+t^2\left(b^2\tilde{\lambda}+b(bpw)'+(b^2pw)'\right)+t\left((pbw)''+b'\tilde{\lambda}\right)\end{equation} $$=t^3b^3h+t^2bk+tl,$$
where $h$, $k$, $l$ are different bounded functions on $[-1,1]$.
Thus the members of \eqref{D} have to be decomposed to three parts and according to \eqref{dif} each part has to be decomposed to further three parts. As we have seen above the computations with the different terms of \eqref{D} are similar, we estimate only the first one, say.
$$\left|D_{n,m}^{\alpha+2,\beta+2, \alpha,\beta,\alpha+2,\beta+2}[L^2(tbpw+\tilde{\lambda}),t]\right|\le \frac{C}{n}\left|I_{(n+1)^*,m}^{\alpha+2,\beta+2,\alpha,\beta,\alpha+2,\beta+2}[L^2(tbpw+\tilde{\lambda}),t]\right|$$ $$+C\left|I_{n,m}^{\alpha+3,\beta+2,\alpha,\beta,\alpha+3,\beta+2}[L^2(tbpw+\tilde{\lambda}),t]\right|=R+M.$$
In view of \eqref{L} we split $R$ and $M$ to three parts. Considering that $u^je^{-u}$ is uniformly bounded and by \eqref{ib}
$$R_1\le \frac{C}{n}\int_{-1}^1e^{-\overline{Q}t}b^3t^3\left|p_{(n+1)^*}^{\alpha+2,\beta+2}p_{m}^{\alpha,\beta}\right|w^{\alpha+2,\beta+2}dx$$
$$\le \frac{C}{n}\int_{-1}^1\left|p_{(n+1)^*}^{\alpha+2,\beta+2}\right|w^{\frac{\alpha}{2}-\frac{5}{4},\frac{\beta}{2}+\frac{7}{4}}dx\le C n^2.$$
Similarly
$$R_2\le \frac{C}{n}\int_{-1}^1e^{-\overline{Q}t}bt^2\left|p_{(n+1)^*}^{\alpha+2,\beta+2}p_{m}^{\alpha,\beta}\right|w^{\alpha+2,\beta+2}dx$$
$$\le \frac{C}{n}\int_{-1}^1\left|p_{(n+1)^*}^{\alpha+2,\beta+2}p_{m}^{\alpha,\beta}\right|w^{\alpha-1,\beta+2}dx\le C n^2.$$
Finally
$$R_3\le \frac{C}{n}\int_{-1}^1\left|p_{(n+1)^*}^{\alpha+2,\beta+2}p_{m}^{\alpha,\beta}\right|w^{\alpha,\beta+2}dx\le C .$$
The same to the $M$ terms:
$$M_1\le C \int_{-1}^1e^{-\overline{Q}t}b^3t^3\left|p_n^{\alpha+3,\beta+2}p_{m}^{\alpha,\beta}\right|w^{\alpha+3,\beta+2}dx$$
$$\le \int_{-1}^1\left|p_n^{\alpha+3,\beta+2}\right|w^{\frac{\alpha}{2}-\frac{1}{4},\frac{\beta}{2}+\frac{7}{4}}dx\le C n^2,$$
$$M_2\le \int_{-1}^1e^{-\overline{Q}t}bt^2\left|p_{n}^{\alpha+3,\beta+2}p_{m}^{\alpha,\beta}\right|w^{\alpha+3,\beta+2}dx$$
$$\le \int_{-1}^1\left|p_{n}^{\alpha+3,\beta+2}p_{m}^{\alpha,\beta}\right|w^{\alpha,\beta+2}dx\le C n^2,$$
And
$$M_3\le \int_{-1}^1\left|p_{n}^{\alpha+3,\beta+2}p_{m}^{\alpha,\beta}\right|w^{\alpha+1,\beta+2}dx\le C .$$
As it is shown the terms with $[2xLf,t]$ are less than the previous ones. Thus by \eqref{k2} and \eqref{D}
\begin{equation}\label{K2}|K_2(n+1,m)-K_2(n,m)|\le \frac{C}{|n-m|^2},\end{equation}
and the proof of \eqref{l3} is finished.

\proof (of Theorem \ref{t1})
Since
$$W_*^{\alpha,\beta,e}f(n)\le \left\|\sum_{m=0 \atop m\neq n}^\infty f(m)K_t^{\alpha,\beta,e}(n,m)\right\|_{\infty, [0,\infty)}+ \|f(n)K_t^{\alpha,\beta,e}(n,n)\|_{\infty, [0,\infty)},$$
for the first term we can apply \eqref{l1} and \eqref{l3} of Lemma \ref{le1} and Theorem A and for the second one \eqref{l2} of Lemma \ref{le1}.

\section{Discrete heat semigroup associated with Dunkl-Jacobi polynomials}
\subsection{Dunkl-Jacobi polynomials}

Dunkl-Jacobi operators and the corresponding eigenfunctions are defined either on $\mathbb{R}$ or on finite intervals ($(-\pi,\pi)$ or $\left(-\frac{\pi}{2},\frac{\pi}{2}\right)$) are examined by several authors, see e.g. \cite{ro} and the references therein. Harmonic analysis, translation operators, convolution structures are developed using Dunkl-Jacobi operators. Here we give the associated heat semigroup and norm estimates for the maximal operator as above.

Let
$$\varphi_k(t):=\frac{P_k^{\alpha,\beta}(\cos t)}{P_k^{\alpha,\beta}(1)}=\frac{\varrho_k^{\alpha,\beta}p_k^{\alpha,\beta}(\cos t)}{P_k^{\alpha,\beta}(1)},$$
where $P_k^{\alpha,\beta}(1)=\binom{k+\alpha}{k}$. The weight function on $[-\pi, \pi]$ is
$$A_{\alpha,\beta}(t)=(1- \cos t)^\alpha (1+\cos t)^\beta |\sin t|.$$
Considering the operator (acting on a function $f$)
$$L_{\alpha,\beta}[f]:=\frac{1}{A_{\alpha,\beta}}\left(A_{\alpha,\beta}f'\right)', $$
$\varphi_k(t)$ are the eigenfunctions the operator $L_{\alpha,\beta}$ and fulfil the initial-value problem
$$L_{\alpha,\beta}[\varphi_k]=-\lambda_k^2\varphi_k, \ws \ws \ws \varphi_k(0)=1,$$
where
$$\lambda_k:=\lambda_k^{\alpha,\beta}=\sqrt{k(k+\varrho)}.$$
Let $k\ge 0$ and define
$$e_k:=\varphi_k+\frac{1}{i\lambda_k}\varphi_k',\ws\ws e_{-k}:=\overline{e}_k.$$
Consider the operator (acting on $f$)
$$\Lambda_{\alpha,\beta}[f]:=f'+\frac{A_{\alpha,\beta}'}{A_{\alpha,\beta}}f_o,$$
where
$$f_o=\frac{f(t)-f(-t)}{2}$$
is the odd part of the function $f$. Then
\begin{equation}\label{djo}\Lambda_{\alpha,\beta}[e_k]=i\ws\mathrm{sign}\ws k \ws\lambda_k e_k, \ws \ws \ws e_k(0)=1.\end{equation}
Moreover the system $\{e_k\}$ is complete and orthogonal on $[-\pi, \pi]$ with respect to $A_{\alpha,\beta}$, see \cite{v} and the references therein.

Similarly to \eqref{si} $L_{\alpha,\beta}$ can be expressed by $\Lambda_{\alpha,\beta}$:
\begin{equation}\label{DJ} L_{\alpha,\beta}[f]=\Lambda_{\alpha,\beta}^2[f]- \left(\frac{A_{\alpha,\beta}'}{A_{\alpha,\beta}}\right)'f_o.\end{equation}

Let us denote by $\psi_k$ the orthonormalized system, that is if $k\ge 0$ (with $p_{-1}^{\alpha,\beta}\equiv 0$)
$$\psi_k(t):=\psi_k^{\alpha,\beta}(t)=\frac{e_k(t)}{\|e_k\|_{2,A_{\alpha,\beta}}}=\frac{1}{\sqrt{2}}\left(p_k^{\alpha,\beta}(\cos t)+ip_{k-1}^{\alpha+1,\beta+1}(\cos t)\sin t\right),$$
and
$$\psi_{-k}(t)=\overline{\psi_k(t)}.$$
The orthonormal Jacobi polynomials fulfil the three-term recurrence relation
$$xp_n^{\alpha,\beta}(x)=a_{n+1}p_{n+1}^{\alpha,\beta}(x)+b_np_n^{\alpha,\beta}(x)+a_np_{n-1}^{\alpha,\beta}(x),$$
where (recalling the notation $\varrho=\alpha+\beta+1$)
$$a_n=a_n^{\alpha,\beta}=\frac{2}{2n+\varrho-1}\sqrt{\frac{n(n+\alpha)(n+\beta)(n+\varrho-1)}{(2n+\varrho)(2n+\varrho-2)}},$$ $$ b_n=b_n^{\alpha,\beta}=\frac{\beta^2-\alpha^2}{(2n+\varrho+1)(2n+\varrho-1)},$$
see \cite[(4.5.1)]{sz} and subsequently we use the abbreviation
$$A_n:=a_n^{\alpha+1,\beta+1}, \ws \ws \ws B_n:=b_n^{\alpha+1,\beta+1}.$$
Consequently it can be readily derived  that $\psi_k$ fulfils the six-term formula
\begin{equation}\label{recpsi}\cos t \psi_{\pm k}(t)=\frac{1}{2}(a_k\pm A_{k-1})\psi_{k-1}(t)+\frac{1}{2}(a_k\mp A_{k-1})\psi_{-(k-1)}(t)+\frac{1}{2}(b_{k}\pm  B_{k-1})\psi_{k}(t)$$ $$+\frac{1}{2}(b_{k}\mp B_{k-1})\psi_{-k}(t)+\frac{1}{2}(a_{k+1}\pm A_{k})\psi_{k+1}(t)+\frac{1}{2}(a_{k+1}\mp A_{k})\psi_{-(k+1)}(t),\end{equation}
where $k\in \mathbb{N}$, $a_{-1},\dots =0$.

With a different normalization (and with some misprints) the formula above is given in \cite[Theorem 3.3]{c}.\\
Rearranging the orthonormal system $\{\psi_k\}_{k\in\mathbb{Z}}$ as $\psi_0, \psi_1, \psi_{-1}, \psi_{2}, \psi_{-2}, \dots$, we get the matrix of the multiplication operator $M_{\cos t}$ in three-block-diagonal form. For simplicity let us denote by $c_k=\frac{a_k+A_{k-1}}{2}$, $c_k^*=\frac{a_k-A_{k-1}}{2}$ and by $d_k=\frac{b_k+B_{k-1}}{2}$, $d_k^*=\frac{b_k-B_{k-1}}{2}$.
\begin{equation}\label{mt}M_{\cos t}=\left[\begin{array}{cccccccccccc}d_{0}&d_{0}^*&c_1&c_1^*&0&0&0&0&\dots&\dots&0&0\\
d_0^*&d_0&c_1*&c_1&0&0&0&\dots\\ \vdots&\vdots&\dots&\ddots&\vdots&\dots&\vdots&\dots&\dots&0&0&\dots\\
\dots &0&0&c_k&c_k^*&d_k&d_k^*&c_{k+1}&c_{k+1}^*&0&0&\dots\\
\dots &0&0&c_k^*&c_k&d_k^*&d_k&c_{k+1}^*&c_{k+1}&0&0&\dots\\
\dots&\dots&\dots &0&0&c_{k+1}&c_{k+1}^*&d_{k+1}&d_{k+1}^*&c_{k+2}&c_{k+2}^*&0\\
\dots&\dots&\dots &0&0&c_{k+1}^*&c_{k+1}&d_{k+1}^*&d_{k+1}&c_{k+2}^*&c_{k+2}&0\\
\dots& \dots&\dots& \dots&\vdots&0&\dots& \dots& \dots&\vdots&\dots&\dots\end{array}\right].\end{equation}

\subsection{The discrete heat semigroup and the corresponding maximal operator}
Let
\begin{equation}K^{\alpha,\beta,D}_t(n,m):= \int_{-\pi}^\pi e^{-(1-\cos \tau)t} \psi_n(\tau)\overline{\psi_m}(\tau)A_{\alpha,\beta}(\tau)d\tau. \end{equation}
Note that the imaginary parts are odd, thus similar to the recurrence coefficients, the kernel function is also real (and symmetric).

The corresponding operator for an $f\in l^2(\mathbb{Z})$ is defined
\begin{equation}W_t^{\alpha,\beta,D}f(n):=\sum_{m\in\mathbb{Z}}f(m)K^{\alpha,\beta,D}_t(n,m).\end{equation}
As in the previous section, the operator is not positivity preserving.

Taking into consideration \eqref{mt} and the fact that $c_k$ tends to $\frac{1}{2}$ and $c_k^*$, $d_k$ $d_k^*$ tend to zero,
 $\{W_t^{\alpha,\beta,D}\}_{t\ge 0}$ is the discrete heat semigroup (in the sense mentioned above).\\

As previously we can extend the maximal operator as follows.

\begin{theorem}\label{t2} Let $\alpha, \beta\ge-\frac{1}{2}$. The maximal operator of the discrete heat semigroup associated with Dunkl-Jacobi polynomials, $W_*^{\alpha,\beta,D}$, fulfils that\\
\rm{(1)} if $1<p<\infty$ and $w\in A_p(\mathbb{Z})$, then for all $f\in l^2(\mathbb{Z})\cap l^p(\mathbb{Z},w)$
$$\|W^{\alpha,\beta,D}_*f\|_{p,w}\le C\|f\|_{p,w},$$
where $C$ is a constant independent of $f$. That is the operator $W^{\alpha,\beta,D}_*$ can be extended uniquely to a bounded operator from $l^p(\mathbb{Z},w)$ into itself.\\
\rm{(2)} If $w\in A_p(\mathbb{Z})$, then for all $f\in l^2(\mathbb{Z})\cap l^1(\mathbb{Z},w)$
$$\|W^{\alpha,\beta,D}_*f\|_{(1,\infty),w}\le C\|f\|_{1,w},$$
where $C$ is a constant independent of $f$. That is the operator $W^{\alpha,\beta,D}_*$ can be extended uniquely to a bounded operator from $l^1(\mathbb{Z},w)$ into $l^{1,\infty}(\mathbb{Z},w)$.
\end{theorem}

\medskip

As above, it is enough to prove the next lemma.

\begin{lemma} Let $\alpha,\beta \ge-\frac{1}{2}$, $n,m \in\mathbb{Z}$. Then
\begin{equation}\label{ld2}|K^{\alpha,\beta,D}_t(n,n)|\le C,\end{equation}
\begin{equation}\label{ld1}|K^{\alpha,\beta,D}_t(n,m)| \le \frac{C}{|n-m|}, \ws \ws \mbox{if} \ws \ws n\neq m ,\end{equation}
\begin{equation}\label{ld3}|K^{\alpha,\beta,D}_t(n+1,m)-K^{\alpha,\beta,e}_t(n,m)| \le \frac{C}{|n-m|^2}, \ws \ws \mbox{if} \ws \ws n\neq m,m\pm 1,\ws \frac{m}{2}<n<2m, \end{equation}
where $C$ is a constant which may be different at each occurrence, and is independent of $t$, $n$ and $m$.
\end{lemma}

 \proof
 Let us introduce the notation
 $$K^{\alpha,\beta,J}_t(n,m)=\int_{-1}^1e^{-(1-x)t}p_n^{\alpha,\beta}(x)p_m^{\alpha,\beta}(x)w^{\alpha,\beta}(x)dx,$$
 the kernel with respect to the Jacobi polynomials. By the substitution $x=\cos t$,
\begin{equation}K^{\alpha,\beta,D}_t(n,m)=K^{\alpha,\beta,J}_t(n,m)+\mathrm{sign}\ n \ \mathrm{sign} \ m \ K^{\alpha+1,\beta+1,J}_t(n-1,m-1).\end{equation}
Thus \cite[(22), (23), (24)]{acl} imply the result.

\medskip

\proof (of Theorem \ref{t2}) We can proceed as in proof of Theorem \ref{t1}.

\medskip

\noindent \small{Department of Analysis, \newline
Budapest University of Technology and Economics}\newline

\small{ g.horvath.agota@renyi.hu}
\end{document}